\documentclass{amsart}
\usepackage[pagebackref]{hyperref}

\title{Small Norms in Quadratic Fields}
\author{Franz Lemmermeyer}
\address{M\"orikeweg 1, 73489 Jagstzell, Germany}
\email{hb3@ix.urz.uni-heidelberg.de}

\newtheorem{lemma}{Lemma}
\newtheorem{prop}[lemma]{Proposition}

\newcommand{\Z}{{\mathbb Z}}
\newcommand{\Q}{{\mathbb Q}}

\newcommand{\OO}{{\mathcal O}}
\newcommand{\fp}{\mathfrak p}

\begin{document}

\maketitle

\section{Introduction}

The computation of units in algebraic number fields usually is
a rather hard task. Therefore, families of number fields with
an explicitly given system of independent units have been of
interest to mathematicians in connection with the computation of
\begin{itemize}
\item inhomogeneous minima of number fields;
\item unramified extensions of number fields;
\item capitulation of ideal classes therein;
\item complete solution of Thue and index form equations, etc.
\end{itemize}
The first such number fields have been given by Richaud, namely
the quad\-ratic number fields of R-D-type (named after Richaud and
Degert). These number fields have the form $k =\,\Q(\sqrt{m}\,)$,
where $m=t^2+r$, $|r| \le t$, $r|4t$. In \cite{Ankeny} it was
shown how to prove the class numbers of such fields to be strictly
bigger than 1 by making use of a lemma due to Davenport. Hasse
\cite{Hasse} has extended this result to other quadratic fields
of R-D-type without being able to handle all cases. The most
complete effort to this date is by Zhang \cite{Zh2}, who used a continued
fraction approach (see \cite{Zh} for sketches of some of the proofs).
In this paper we intend to show how to prove Zhang's
results using a geometric idea; this has the advantage that it
can be generalized to a family of cubic (and possibly quartic)
fields. For some history and an extensive list of references, see 
the forthcoming book of Mollin \cite{M}.

\section{Quadratic Fields}
In order to show how the proof works in a simple case, we first
will look at quadratic fields $k = \Q(\sqrt{m}\,)$. Let
$\varepsilon>1$ be its fundamental unit; for $\alpha \in$ k  we
let $\alpha '$ denote the conjugate of  $ \alpha $. In particular,
we have $N \xi = \xi \xi'$, where $N = N_{k/\Q}$ denotes the 
absolute norm.

Now let $ \xi \in $ k and a real positive number c be given;
we can find a unit $ \eta \in  E_{k} $ such that
\begin{equation}
c \leq |\,\xi \eta\,|<c \varepsilon.
\end{equation}
Letting $ n=|N\xi|$ we find $ |\xi' \eta'|=n/|\xi \eta|$,
so equation (1) yields
\begin{equation}
\frac{n}{c \varepsilon} \leq |\,\xi' \eta'\,|< \frac{n}{c}.
\end{equation}
Writing $ \alpha = \xi \eta = a+b \sqrt{m} $ gives
$|2a|=|\alpha + \alpha '| \leq |\alpha|+|\alpha '|< c \varepsilon
+ \frac{n}{c},$
and correspondingly,
$|2b\sqrt{m}\,|=|\alpha -\alpha '| \leq |\alpha|+|\alpha'|<
c\varepsilon+\frac{n}{c}.$
Because we want the coefficients $a$ and $b$ to be as small
as possible, we have to choose $c$ in such a way that
$ c\varepsilon + n/c $ becomes a minimum. Putting
$ c = \sqrt{n/\varepsilon} $ we get
\begin{equation}
|2a|< 2\sqrt{n\varepsilon},\quad |2b\sqrt{m}\,|
       < 2\sqrt{n\varepsilon}.
\end{equation}
Making use of a lemma due to Cassels, we can improve these bounds:
\begin{lemma}\label{th1}
Suppose that the positive real numbers $x, y$ satisfy the
inequalities $x \leq s,\ y \leq s$, and $xy \leq t$. Then,
$x+y \leq s+t/s$.
\end{lemma}

\noindent
Proof. $0 \leq (x-s)(y-s) = xy-s(x+y)+s^{2} \leq s^{2}+t-s(x+y)$.

\medskip
Putting $x=|\alpha|$ and $y=|\alpha'|$ in Lemma \ref{th1}
we find
\begin{displaymath}
|2a| \leq |\alpha|+|\alpha'|< \sqrt{n\varepsilon}+
\sqrt{n/\varepsilon},
\end{displaymath}
and likewise
\begin{displaymath}
|2b \sqrt{m}\,|< \sqrt{n\varepsilon}+\sqrt{n/\varepsilon}.
\end{displaymath}
We have proved

\begin{prop}\label{th2}
Let $k = \Q(\sqrt{m}\,)$ be a real quadratic number field,
$\varepsilon > 1$ a unit in k, and $0 \not= n=|N\xi\,|$
for $\xi \in k$. Then there is a unit $\eta=\varepsilon^{j}$
such that $\xi \eta = a+b\sqrt{m}$ and 
\[|a| < \frac{\sqrt{n}}{2}
(\sqrt{\varepsilon}+1/\sqrt{\varepsilon}\,), \qquad
|b| < \frac{\sqrt{n}}{2 \sqrt{m}} (\sqrt{\varepsilon}
+1/\sqrt{\varepsilon}\,).\]
\end{prop}

Suppose that we are looking for a $\xi \in \Z[\sqrt{m}\,]$
with given norm $\pm n$. If we know a unit $\varepsilon >1$,
we can use Proposition \ref{th2} to find a power $\eta$ of
$\varepsilon$ such that $\xi \eta=a+b \sqrt{m}$ has bounded
integral coefficients $a,\ b$. Moreover, the bounds do not depend
on $\xi$. In order to test if a given $n$ is a norm in $k/\Q$,
we therefore have to compute only the norms of a finite number
of elements of $k$. Similar results are valid in case
$\{1,\theta \}$ is an integral basis of the ring $\OO_{k}$ of
integers in $k$, where $\theta = \frac{1}{2}(1+\sqrt{m}\,)$.

After these preparations, it is an easy matter to prove the
following result originally due to Davenport:
\begin{prop}\label{th3}
Let $m, n, t$ be natural numbers such that $m=t^{2}+1\,$; if
the diophantine equation $|x^{2}-my^{2}|=n$ has solutions in
$\Z$ with $n < 2t$, then $n$ is a perfect square.
\end{prop}

\noindent
Proof. Let $\xi = x+y\sqrt{m}$; then $|\,N \xi\,|=n$, and since
$\varepsilon =t+u\sqrt{m}\,>1$ is a unit in $\Z[\sqrt{m}\,]$,
we can find a power $\eta$ of $\varepsilon$ such that
$\xi \eta=a+b\sqrt{m}$ has coefficients $a, \ b$ which satisfy
the bounds in Proposition 2.2. Since $2t<\varepsilon < 2\sqrt{m}$,
we find
\[|b| \leq \frac{\sqrt{n}}{2\sqrt{m}}\, \left( \sqrt{\varepsilon}
  +\frac{1}{\sqrt{\varepsilon}}\right) < 1+\frac{1}{t}.\]
Since the assertion is trivial if $t=1$, we may assume that
$t \ge2$, and now the last inequality gives $|b| \leq 1$.
If $b=0, |N\xi|=a^{2}$ would be a square; therefore,
$b=\pm 1$, and this yields $\alpha = \xi \eta=a \pm \sqrt{m}.$
Now $|N\xi| = |N\alpha| = |a^{2}-m|$ is minimal for
values of $a$ near $\sqrt{m}$, and we find

\begin{tabular}{lccl}
$|a^{2}-m|$ & $=$ & $2t$ & if $a=t-1$;\\
$|a^{2}-m|$ & $=$ & $ 1$ & if $a=t$;  \\
$|a^{2}-m|$ & $=$ & $2t$ & if $a=t+1$.
\end{tabular}

\noindent
This proves the claim.

\medskip
Using the idea in the proof of Proposition \ref{th3} one can
easily show more:

\begin{prop}\label{th4}
Let $m, n, t$  be natural numbers such that $m=t^{2}+1$; if
the diophantine equation $|x^{2}-my^{2}|= n$ has solutions in
$\Z$ with $n < 4t+3$, then $n = 4t-3$, $n = 2t$, or $n$ is a
perfect square.
\end{prop}

In \cite{Ankeny}, Proposition \ref{th3} was used to show that the
ideal class group of $k = \Q(\sqrt{m}\,)$ has non-trivial
elements (i.e. classes that do not belong to the genus class
group) if $m=t^{2}+1$ and $t=2lq$ for $l>1$ and prime $q$:
since $m \equiv 1 \bmod q$, $q$ splits in $k$, i.e. we have
$(q)=\fp\fp''$. If $\fp$ were principal, the equation
$x^{2}-my^{2}=\pm 4q$ would have solutions in $\Z$; but since
$4q < 2t = 4lq$ is no square, this contradicts Proposition \ref{th2}.

If we consider the case $m=t^{2}+2$ instead of $m=t^{2}+1$, then
the method used above does not seem to work: for example, the
equation $x^{2}-my^{2}=-2$ is solvable (put $x=t, y=1$) and 2 is
no square. Therefore, we have to modify our proof in order to get
non-trivial results.
\begin{prop}\label{th5}
Let $m, n, t$ be natural numbers such that $m=t^{2}+2$ and $t \ge 12$;
if the diophantine equation $|x^{2}-my^{2}|=n$ has solutions in $\Z$
and if neither $n$ nor $2n$ are perfect squares, then
$n=2t\pm 1, 4t-7, 4t-2, $ or $n \ge 4t+2$.
\end{prop}

\noindent
Proof.  Let $\xi= x+y\sqrt{m}, n=|N\xi|$,
and suppose that neither $n$ nor $2n$ are perfect squares.
Letting $\delta=t+\sqrt{m}$, we find $\delta^{2}=2\varepsilon$,
where $\varepsilon$ is a unit in $\Z[\sqrt{m}\,]$.
Obviously we can find a power $\eta$ of $\varepsilon$ such that
\[\sqrt{n\sqrt{\varepsilon}}/\varepsilon \le |\xi\eta|
  < \sqrt{n\sqrt{\varepsilon}}.\]
Now we distinguish two cases:
\begin{description}
\item[1.] $\quad \sqrt{n/\sqrt{\varepsilon}} \le |\xi \eta| <
  \sqrt{n \sqrt{\varepsilon}}:$ 
  Writing $\xi\eta = a+b\sqrt{m}$, we find that $|b| \le 1$.
  If $b=0$, then $b$ is a square, so assume $b=\pm 1$.
  Then the same reasoning as in Proposition \ref{th3} shows that
  either $n=2$, i.e. $2n$ is a square, or that $n=2t\pm 1,n=4t-2$
  or $n \ge 4t+2$.
  \item[2.] $\quad \sqrt{n\sqrt{\varepsilon}}/\varepsilon
   \le |\,\xi \eta\,| < \sqrt{n\sqrt{\varepsilon}}:$ 
   Multiplying $\xi \eta$ with $\delta$ we get
   \[\sqrt{2n/\sqrt{\varepsilon}} \le |\,\xi \eta \delta\,| <
   \sqrt{2n \sqrt{\varepsilon}}.\]
   As in case 1. above, we find $|b| \le 1$, where
   $\xi \eta \delta = a+b\sqrt{m}$; if $b=0$, then
   $2n=|N(\xi \eta \delta)|=a^{2}$ is a square. If $b=\pm 1$, then
   $a^{2}-mb^{2}$ must be even (because $|N\xi\delta|$ is even), 
   or $|N\xi\delta| = 4t \pm 2, =8t-14$, or
   $\ge 8t+14$ (because $12t-34 \ge 8t+14$ for all $t \ge 12$).
\end{description}

\medskip
Exactly as after the proof of Proposition \ref{th3}, we can deduce
results about the ideal class group of $\Q(\sqrt{m}\,)$ from
Proposition \ref{th5} if $m=t^{2}+2$. Moreover we remark that we
have treated these two cases only to explain the method; similar
results can be proved for other quadratic fields of R-D-type.
As an example, we give the corresponding result for $m = t^2-2$:

\begin{prop}\label{tz}
Let $m, n, t$ be natural numbers such that $m=t^{2}-2$ and $t \ge 12$;
if the diophantine equation $N(\xi) = |x^{2}-my^{2}|=n$ has solutions
$\xi = a+b\sqrt{m} \in \Z[\sqrt{m}\,]$, then either $\xi = n\eta$
for some $n \in \Z$ and some unit $\eta \in Z[\sqrt{m}\,]$, or
$n=2t\pm 3, 4t-9, 4t \pm 6, $ or $n \ge 4t+6$, and $\xi$ is associated
to  one of $\{t \pm 1 \pm \sqrt{m}, \, t \pm 2 \pm \sqrt{m},\,
2t - 1 \pm 2\sqrt{m},\, 2t \pm 2 \pm 2\sqrt{m} \}$.
\end{prop}


\begin{thebibliography}{ACH}

\bibitem[ACH]{Ankeny}
N. C. Ankeny, S. Chowla, H. Hasse, {\em On the class number of
the real subfield of a cyclotomic field}, J. Reine Angew. Math.
{\bf 217} (1965), 217--220
%

\bibitem[H]{Hasse}
H. Hasse, {\em \"{U}ber mehrklassige, aber eingeschlechtige
reell-quadratische Zahlk\"{o}rper}, El. Math. {\bf 20} (1965), 49--59
%


\bibitem[M]{M} R. Mollin,
{\em Quadratics}, CRC
%


\bibitem[Zh]{Zh}
Xian-ke Zhang, {\em Determination of solutions and solvabilities of
diophantine equations and quadratic fields}, preprint
%

\bibitem[Zh2]{Zh2}
Xian-ke Zhang, {\em Solutions of the diophantine equations related
to real quadratic fields}, Chin. Sci. Bull {\bf 37} (1992), 885--889
%

\end{thebibliography}
\end{document}